%% file: asymptoticfullerene.tex
 \documentclass[11pt]{article}
 \usepackage[top=1in,bottom=1in,left=1.5in,right=1.5in]{geometry}
  \usepackage{amsmath,amssymb,amsthm}
  \usepackage{graphicx}
  \usepackage{epsfig}
  \usepackage{amsmath,amssymb}
  \usepackage{subcaption}
  \usepackage{color}
  \usepackage{hyperref}
  \newcommand{\hs}{}
  \newtheorem{theo}{\hs Theorem}[section]

  \newtheorem{rem}[theo]{\bfseries \hs Remark}
  
  \numberwithin{equation}{section} 

\begin{document}

\title{Asymptotic enumeration of perfect matchings in $m$-barrel fullerene graphs}

\author{%
   Afshin Behmaram\thanks{Faculty of Mathematical Sciences, University of
   Tabriz, Tabriz, Iran.\newline \emph{email}: behmaram@tabrizu.ac.ir},
C{\'e}dric Boutillier\thanks{Laboratoire de Probabilit\'es et Mod\`eles
  Al\'eatoires, UPMC Univ. Paris 06, 4 place
Jussieu, F-75005 Paris, France.\newline \emph{email}: cedric.boutillier@upmc.fr}
}

\date{\today}
\maketitle

\begin{abstract}
A connected planar cubic graph is called an $m$-barrel fullerene and denoted by
$F(m,k)$, if it has the  following structure:
The first circle is an $m$-gon. Then $m$-gon  is bounded by $m$ pentagons.
After that  we have additional  k layers of hexagons. At the last circle
$m$-pentagons connected to the second $m$-gon. In this paper we asymptotically
count by two different methods the number of perfect matchings in $m$-barrel
fullerene graphs, as the number of hexagonal layers is large, and
show that the results are equal.
\end{abstract}

\noindent {\bf 2010 Mathematics Subject Classification.} 05C30, 05C70, 15A15.

\noindent {\bf Keywords.}  perfect matchings,  fullerene graph, m-barrel fullerene

\section{Introduction}

A \emph{fullerene graph} is a cubic, planar, 3-connected graph with only
pentagonal and hexagonal faces. It follows easily from the Euler's formula
that there must be exactly 12 pentagonal faces, while the number of hexagonal
faces can be zero or any natural number greater than one. The smallest
possible fullerene graph is the dodecahedron on 20 vertices, while the
existence of fullerene graphs on an even number of vertices greater than 22
follows from a result by Gr\"unbaum and Motzkin \cite{grunbaum}. Classical
fullerene graphs have been intensely researched since the discovery of
buckminsterfullerene in the fundamental paper \cite{kroto}, which appeared
in 1985, and gave rise to the whole new area of fullerene science.

A connected $3$-regular planar graph $G=(V,E)$ is called an
\emph{$m$-generalized fullerene} if exactly two of its faces are $m$-gons and all
other faces are pentagons and/or hexagons. (We also count the outer (unbounded)
face of $G$.) In the rest of the paper we only consider $m \geq 3$;
note that for $m=5,6$ an $m$-generalized fullerene graph is
a classical fullerene graph. As for the classical fullerenes it is easy to show
that the number of pentagons is fixed, while the number of hexagons is not
determined. The smallest  $m$-generalized fullerene has $4m$ vertices and
no hexagonal faces. By inserting
$k \geq 0$ layers of $m$ hexagons between two layers of pentagons we reach the
symmetric class of m-generalized fullerenes called \emph{$m$-barrel
fullerenes}.

\begin{figure}[ht]
    \centering
    \includegraphics[width=5cm]{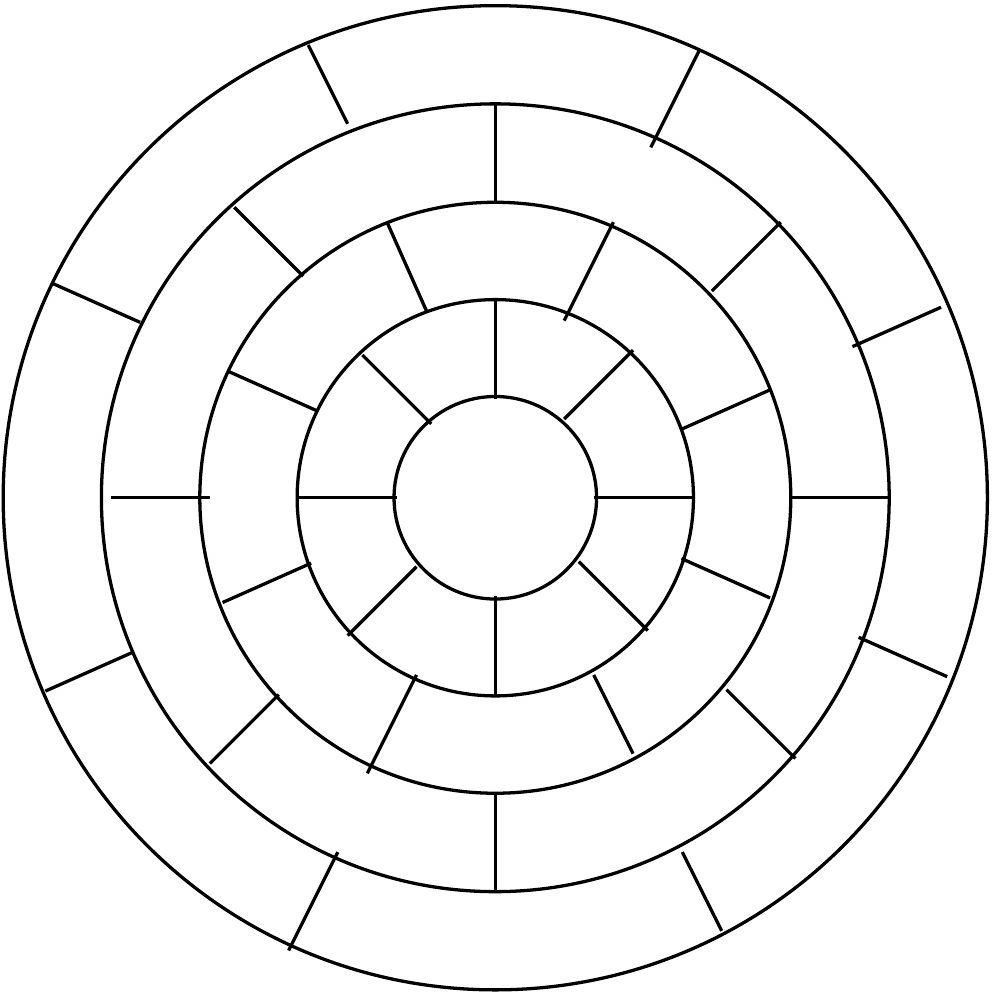}
\caption{The $m$-barrel fullerene $F(8,2)$.}
\label{circles}
\end{figure}

The $m$-barrel
fullerene with $k$ layers of hexagons, denoted by $F(m,k)$, can be defined as a
sequence of concentric layers as follows: the first circle is an
$m$-gon. This $m$-gon  is bounded by $m$ pentagons.  After that  we have
additional $k$ layers of $m$ of hexagon. Then one again has a circular layer
with  $m$-pentagons  connected to the second $m$-gon, represented by the outer
face.
$m$-barrell fullerenes can be neatly represented graphically using a sequence of
$k+3$ concentric circles with monotonically increasing radii
such that the innermost and the outermost circle
each have $m$ vertices (representing, hence, two $m$-gons), while all other
circles have $2m$ vertices each, connecting alternatively to vertices of the
larger or smaller circle to create hexagonal an pentagonal faces. An example is
shown in Figure~\ref{circles}.

The $m$-barrel fullerenes are the main subjects of the present paper, since
their highly symmetric structure allows for obtaining good bounds and even exact
results
on their quantitative graph properties. For example Kutnar and
Maru{\v s}i{\v c} in \cite{KM08} studied Hamiltonicity and cyclic
edge-conectivity of $F(5,k)$. See also~\cite{BF1} for some structural results
about $m$-barrel fullerene graphs, such as the diameter, Hamiltonicity and the
leapfrog transformation.

A \emph{matching} $M$ in a graph $G$ is a collection of edges of $G$ such that
no two edges of $M$ share a vertex. If every vertex of $G$ is incident to an
edge of $M$, the matching $M$ is said to be \emph{perfect}.
A perfect matching is also often called a \emph{dimer configuration} in mathematical
physics and chemestry.
Perfect matchings have
played an important role in the chemical graph theory, in particular for
benzenoid graphs, where their number correlates with the compound's stability.
Although it turned out that for fullerenes they do not have
the same role as for benzenoids, there are many results concerning their
structural and enumerative properties. See
\cite{BF12,andova1,dosliclff,doslicfmm,KKMS09} for more result
on perfect matchings in fullerenes. We denote by $\Phi(G)$ the number of perfect
matchings of $G$.

The goal of this paper is to compute the growth constant $\rho(m)$ for the number of
perfect matchings for the family of graphs $F(m,k)$ for a fixed $m$, as $k$ goes
to infinity:
\begin{equation}
  \rho(m) = \lim_{k\to\infty} \Phi(F(m,k))^{1/k}.
  \label{eq:growth}
\end{equation}

Behmaram, Doslic and Friedland in \cite{BF1} obtained some exact results for
the number of perfect matching for the values of $m\leq 5$, using the transfer
matrix method, described in Section~\ref{sec:transfer}, allowing in particular
for a direct computation of the growth $\rho(m)$ for $m\in\{3,4,5\}$, see
Theorem~\ref{thm:BF1}.

In this paper, we estimate the number of perfect
matchings of $F(m,k)$ for large $k$, and compute exactly $\rho(m)$ for any $m$.

\begin{theo}
  Let $m\geq 3$.   The growth constant for the family of $m$-barrel fullerenes
  is equal to
  \begin{equation*}
    \rho(m)=
    \prod_{j=1}^{\lfloor\frac{m+1}{3}\rfloor}\left(2\cos
    \frac{\pi(2j-1)}{m}\right)^2.
  \end{equation*}
\end{theo}

We propose two proofs using two different methods from combinatorics:
\begin{itemize}
  \item the transfer matrix approach (Section~\ref{sec:transfer}), which is
    explicitly diagonalized using (a baby version) of the Bethe
    Ansatz~\cite{Bethe1931}, method from integrable systems to diagonalize the
    Hamiltonian of integrable spin chains for example;
  \item an approach using coding with non-intersecting paths, counted by
    determinants via the Lindstr\"om-Gessel-Viennot lemma.
\end{itemize}

These two approaches are classical in some branches of combinatorics, and
believe that they can be of great use to study properties of fullerene graphs,
especially those with some symmetry.

Before we introduce these methods, we apply some transformations to the graph
$F(m,k)$ to see it as piece of the hexagonal lattice wrapped on the cyclinder,
with specific boundaray conditions, in order to finally reformulate the question
of perfect matchings on $F(m,k)$ as a problem of tilings with rhombi.

\section{Perfect matchings on $m$-barrel fullerenes and tilings of cylinders
with rhombi}
\label{sec:tiling}


\begin{figure}[ht!]
  \centering
  \begin{subfigure}[c]{4cm}
    \includegraphics[width=4cm]{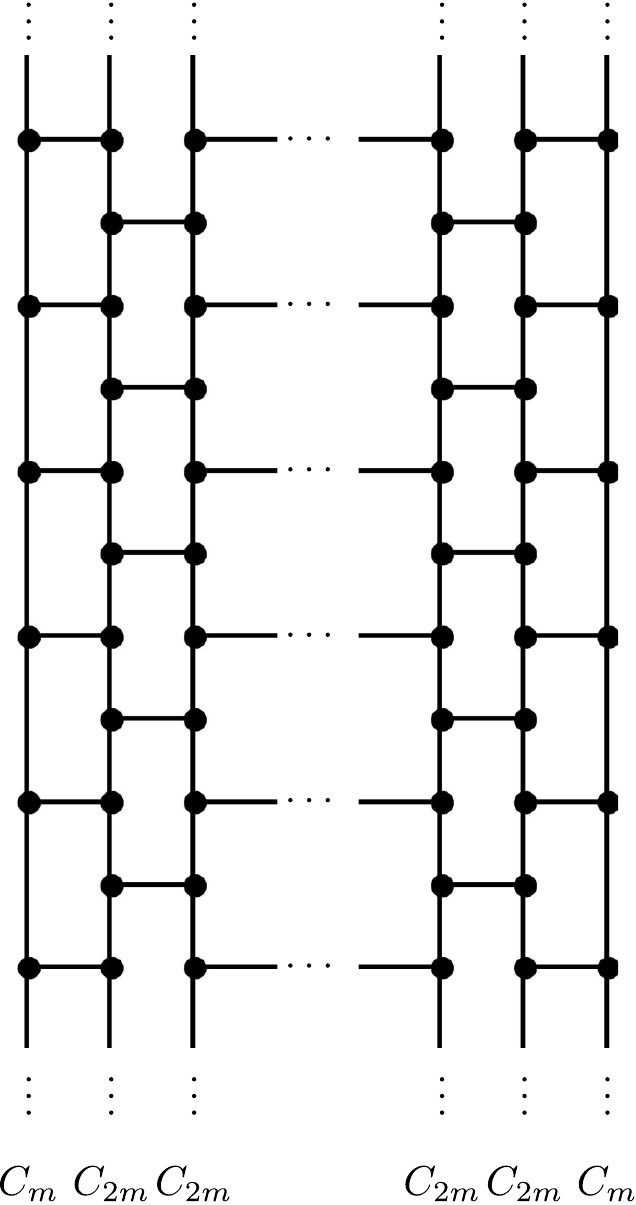}
  \end{subfigure}
  \hspace{0.2\textwidth}
  \begin{subfigure}[c]{4cm}
    \includegraphics[width=4cm]{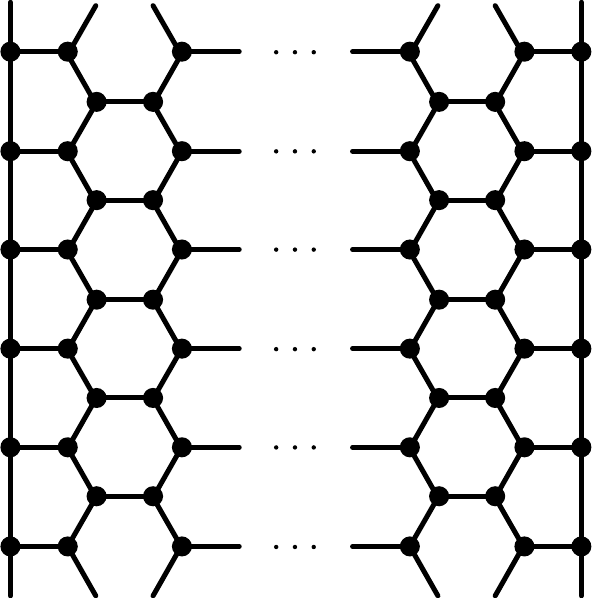}
  \end{subfigure}
  \caption{Left: an $m$-barrel fullerene, represented on the cylinder. The top and
    bottom are identified to create the cylinder. The cycles $C_m$ on the
    extremities are the boundary of the cylinder and correspond to the boundaries of
    the $m$-gons in the representation of Figure~\ref{circles}. Right: the same
  graph on the cylinder, slightly deformed so that hexagonal faces are regular.
  \label{barrel}}
\end{figure}

To begin with, instead of presenting the graph on the sphere, as on
Figure~\ref{circles}, we draw it on the cylinder, where now the $m$-gons
represent the two components of the boundary of this cylinder. These two cycles
of size $m$ are separated by $k+1$ cycles of size $2m$ winding around the
cylinder, each of them connected to their left and right neighboring cycle by
$m$ horizontal edges to create the pentagonal/hexagonal faces, arranged in a
brickwall pattern on the cylinder. See
Figure~\ref{barrel} (left).
This brickwall pattern can in turn be deformed so that the hexagonal faces are
now regular polygonal, so that the bulk of the graph looks like the regular
honeycomb lattice wrapped on a cylinder. See Figure~\ref{barrel} (right).



\begin{figure}[ht!]
  \centering
  \begin{subfigure}[c]{5cm}
    \includegraphics[width=5cm]{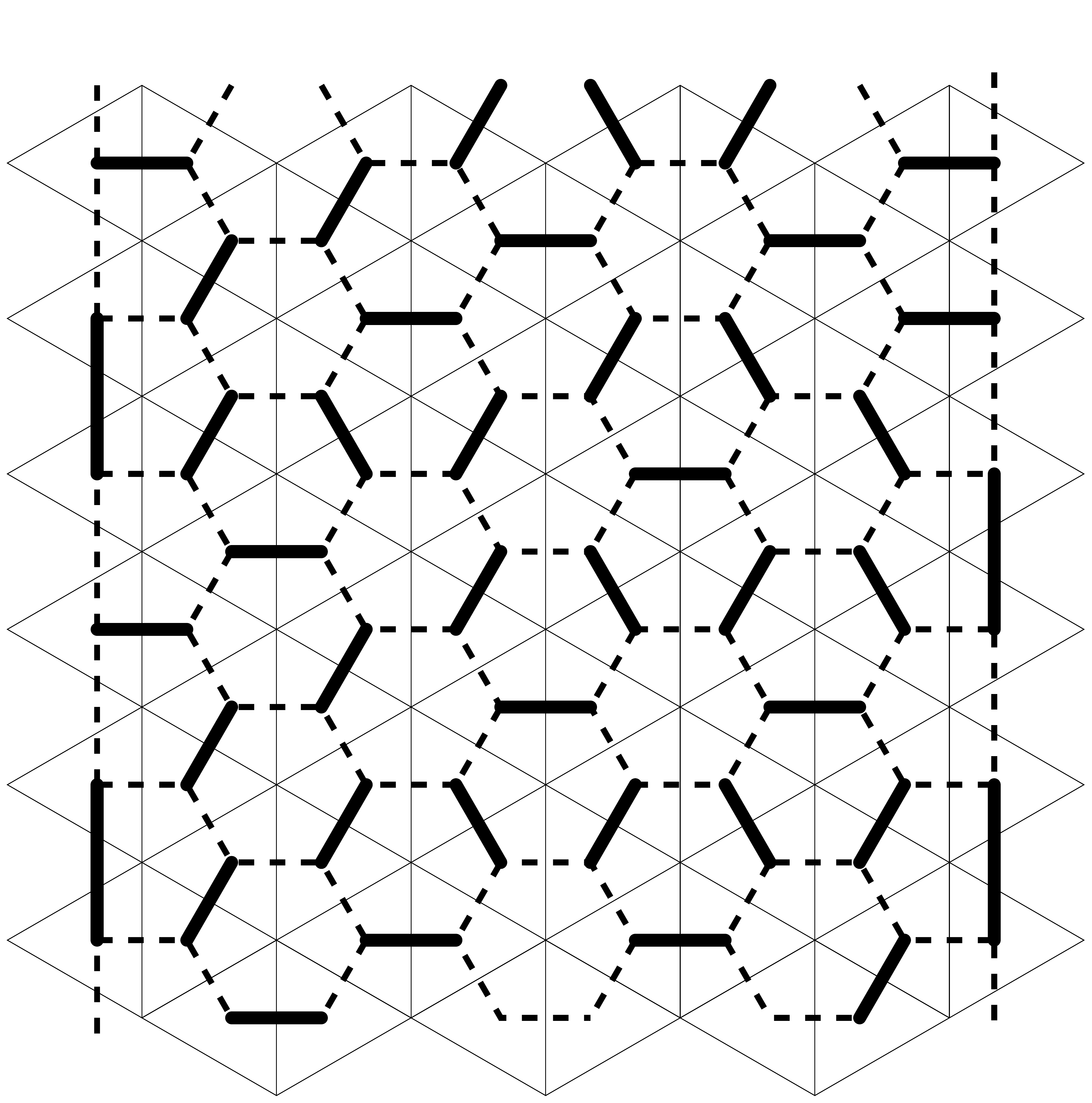}
  \end{subfigure}
  \hspace{0.1\textwidth}
  \begin{subfigure}[c]{5cm}
    \includegraphics[width=5cm]{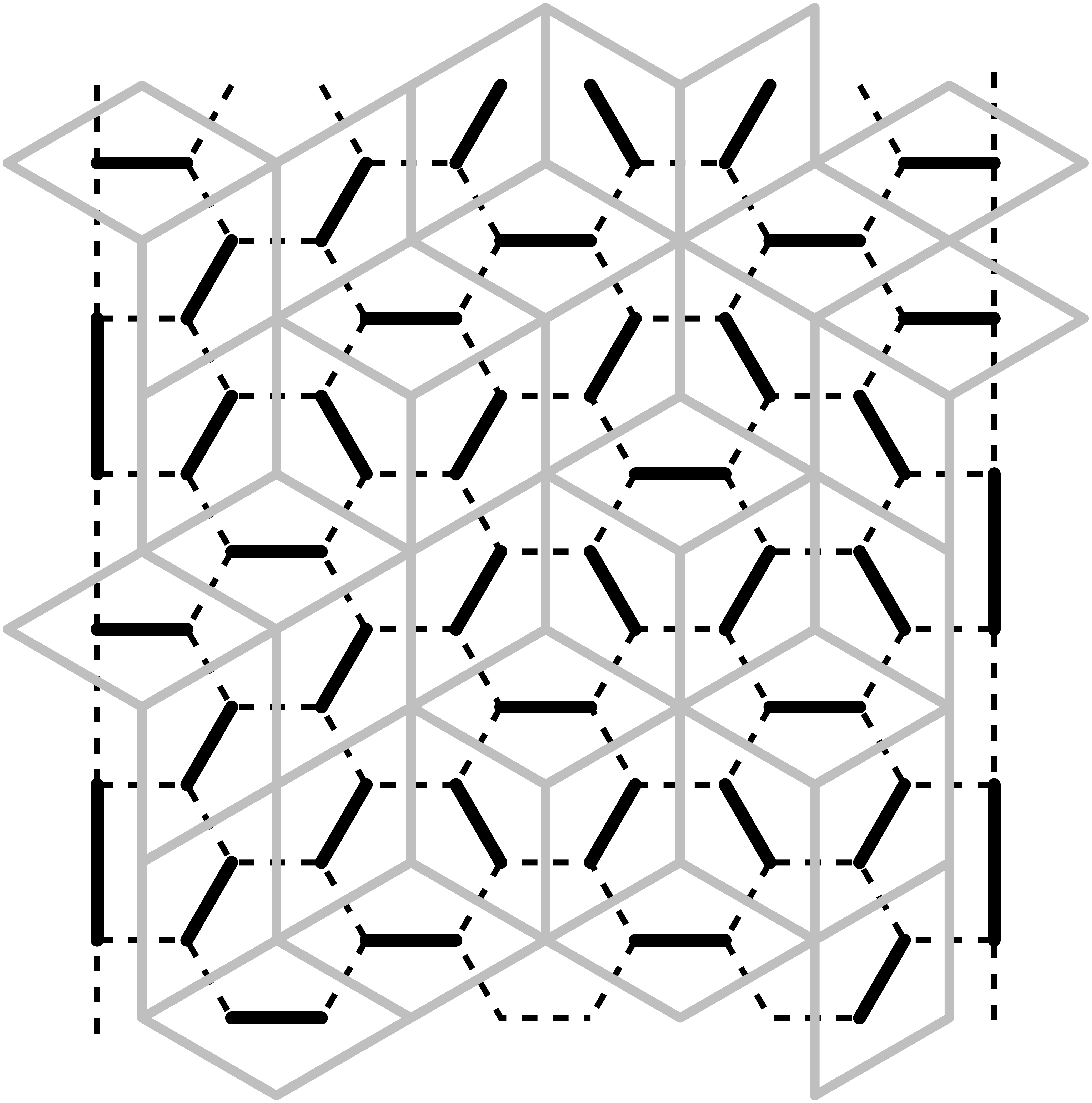}
  \end{subfigure}

    \caption{The correspondence between perfect matchings on $F(6,5)$ (left) and
      rhombi tiling on the cylinder. Recall that the graph is drawn on a cylinder, so
  that the edges at the top are connected to vertices at the bottom.}
    \label{barrel4}
  \end{figure}

It is obvious that perfect matchings such a hexagonal graph on the cylinder with
pentagonal faces on the boundary are in
bijection with tilings of the cylinder with unit rhombi, where some of the
rhombi are allowed to stick out of the boundary. See Figure~\ref{barrel4}.




In the sequel, we will use indifferently the language of perfect matchings and
tilings with rhombi. In particular horizontal edges in a perfect matching of
$F(m,k)$ correspond exactly to horizontal rhombi in the tiling picture. Note
that as long as there is at least a horizontal rhombus in the tiling, the
position of all horizontal rhombi is sufficient to reconstruct the whole tiling.
On the other hand, the two other types of rhombi (with vertical edges) form
non-interesecting paths connecting the left and right boundaries. The collection
of these paths also characterize the whole tiling.

In the next two sections, we present two methods to count the number of perfect
matchings of $F(m,k)$: the transfer matrix method, which uses the horizontal
rhombi, and a method using the collection of non-intersecting paths.

\section{The transfer matrix method}
\label{sec:transfer}

The first method we present, using a \emph{transfer matrix}, is well suited
to count the number of perfect matchings on graphs with
some regularity or periodicity, which is indeed the case here.

Before introducing the transfer matrix, we need to fix some notation.
For $j\in\{1,..,k\}$, denote by $E_j$ the subset of horizontal edges of $F(m,k)$
between the $j$th and the $(j+1)$th cycle $C_{2m}$. Extend this definition to
$E_0$ (resp. $E_{k+1}$) to be the subset of horizontal edges before the first (resp.
after the last) cycle $C_{2m}$, which come from the pentagons at the left (resp.
right) end.

\begin{figure}[hb!]
  \centering
  \resizebox{3cm}{!}{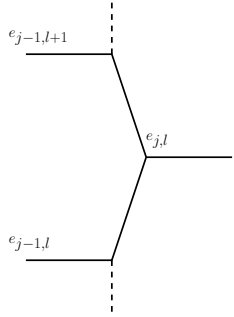}
  \caption{Labeling of the horizontal edges of $F(m,k)$.\label{fig:label_edges}}
\end{figure}

We label the edges of $E_j$ by $e_{j,l}$, $l\in\{0,\dots,m-1\}$. Indices $l$
increase when going ``up'', and they should be thought modulo $m$, so that
$e_{j,l+m}$ is the same as $e_{j,l}$.
For consistency of the labeling between different layers of horizontal edges,
we take the convention that $e_{j,l}$ has on its left $e_{j-1,l}$ just below and
$e_{j-1,l+1}$ just above. See Figure~\ref{fig:label_edges}.
Subsets of $E_j$ are thus in bijection of subsets of $I_m=\{0,\dots,m-1\}$.

Given two subsets $S_{j-1}$ and $S_{j}$ of respectively $E_{j-1}$ and $E_{j}$, let
$a_j(S_{j-1}, S_{j1})$ be the number of perfect matchings of the $j$-th cycle
where we removed vertices attached to the edges in $S_{j-1} \cup S_j$.
Because of invariance by translation, this number does not really depend on $j$,
but only on the subsets $S,T\subset I_m$ in bijection with $S_{j-1}$ and $S_j$
respectively. We store all these numbers in a matrix $A=(a_{S,T})_{S,T\subset
I_m}$. This matrix $A$ (which implicitly depends on $m$) is called the
\emph{transfer matrix} of the model.

It will be convenient to use the \emph{bra} and \emph{ket} notation from quantum
mechanics. The matrix $A$ is thought as the matrix of a linear operator on a
vector space in the orthonormal basis $(|T\rangle)_{T\subset I_m}$ of some
vector space of dimension $2^m$, indexed by subsets of $I_m$.

The dual basis $(\langle S|)_{S\subset I_m}$ satisfies that for all $S, T\subset
I_m$,
\begin{equation*}
  \langle S| T\rangle =
  \begin{cases}
    1 & \text{if $S=T$,} \\
    0 & \text{otherwise.}
  \end{cases}
\end{equation*}

With these notations, we have $A=\sum_{S,T} a_{S,T} |S\rangle\langle T|$ and
$a_{S,T} = \langle S|A| T\rangle$.

Let us make a few remarks on this matrix $A$:
\begin{itemize}
  \item $A$ is invariant by ``vertical translation''.
  \item because of the left/right and top/bottom symmetry, we have $\langle
    S|A|T\rangle=\langle T|A|S\rangle$: in this basis, $A$ is symmetric, thus
    diagonalisable, with orthogonal eigenspaces.
  \item Due to geometric constraints of the problem, two sets $S_{j-1}$ and
    $S_j$ coming from a perfect matching should interlace. In particular if
    $|S|\neq |T|$, $\langle S|A|T\rangle=0$. If we write down the matrix, with
    subsets indexing rows and columns ordered according to their cardinal, then
    $A$ has a block diagonal structure.
  \item Perron-Frobenius theorem guarantees, that on each block, the largest
    eigenvalue is non-degerate, and associated with an eigenvector wich can be
    chosen with all its entries in that block to be positive.
\end{itemize}

In this vector space, one can also define a \emph{boundary vector}
$|\Omega\rangle$,
corresponding to a formal linear combination of possible configurations coming
from possible perfect matchings of $F(m,k)$. Then it is a simple observation
that the number of perfect matchings $\Phi(F(m,k))$ of $F(m,k)$ can be expressed
with $A$ and $|\Omega\rangle$ as follows
\begin{equation}
  \Phi(F(m,k)) = \langle \Omega | A^{k+1} | \Omega \rangle,
  \label{eq:transfer_mat}
\end{equation}
as there are $k+1$ transitions between $E_{j-1}$ and $E_j$, for $1\leq j \leq
k+1$, each of of them being encoded by a matrix $A$.

Computing explicitly the reduced form of $A$ for $m=3,4$ and $5$, the authors
of~\cite{BF1} gave an exact formula for
$\Phi(F(3,k)$,
$\Phi(F(4,k)$, and
$\Phi(F(5,k)$:

\begin{theo}[\cite{BF1}]
  \label{thm:BF1}
let $\Phi (F(m,k))$ denote the number of perfect matching in $F(m,k)$ then for
m=3,4,5 we have:
\begin{align*}
  \Phi (F(3,k)) &= 3^{k+2}+1, \\
  \Phi (F(4,k))  &= 2(2+\sqrt{2})^{k+1} + 2(2-\sqrt{2})^{k+1} + 2^{k+3} + 1,\\
  \Phi (F(5,k)) &= 5^{k+2} + 5 \left [ \left ( \frac{5 + \sqrt{5}}{2} \right )^{k} + \left ( \frac{5 - \sqrt{5}}{2} \right )^{k} \right ] + 1.
\end{align*}
\end{theo}

\subsection{The Bethe Ansatz}

It turns out that it is possible to express the eigenvalues (and eigenvectors)
of $A$ for any value of $m$, using the so-called \emph{Bethe
Ansatz}~\cite{Bethe1931}, a method
used in theoretical physics to diagonalize the Hamiltonian of integrable systems
with interaction, by looking for eigenvectors as a superposition of plane waves.

In our problem, the situation is particularly simple: it turns out, as we will
see later, that all the eigenvectors are expressed in terms of determinants
(Slater determinants in quantum mechanics terminology).  This is a feature of
the dimer model and its free-fermionic nature.

As mentioned above, $A$ has a block diagonal structure. The blocks (or
\emph{sectors}) are indexed
by $p\in\{0,\dots,m\}$, the number of elements of the corresponding subsets
indexing rows and columns. The block indexed with $p$ has size $\binom{m}{p}$.

For any $p$, we look for $\binom{m}{p}$ eigenvectors as linear combination of
basis vectors $|S\rangle$, where $|S|=p$.

For the counting problem we are interested in, the structure of eigenvectors and
eigenvalues is a bit degenerated. In order to perform the computations, it is
easier to introduce a weighted version of the transition matrix:
Fix $b$, and $c$ positive, distinct real numbers, and put weight $b$ (resp. $c$)
on every upgoing (downgoing) edge from left to right.
let $B=(\langle S|B|T\rangle)$ be the
matrix such that $\langle S|B|T\rangle$ is the sum of all perfect matchings
on a layer with boundary conditions described by $S$ and $T$: if $S$ and $T$ are
empty, then $\langle S|B|T\rangle = b^m + c^m$. If $S$ and $T$ are not empty
and compatible, then there is only one perfect matching for the transition
from $S$ to $T$ and $\langle S|B|T\rangle$ is of the form $b^j c^{m-j-|S|}$ for
some $j$.

\subsection{The $p=0$ and $p=1$ sectors}

The sector $p=0$ is one-dimensional, and spanned by $\emptyset$. If we don't
remove any vertices, the cycle of length $2m$ has 2 perfect matchings,
consisting of odd and even edges respectively. Therefore, we have
\begin{equation*}
  B|\emptyset\rangle = (b^m+c^m) |\emptyset\rangle.
\end{equation*}

The sector $p=1$ is $m$-dimensional. For $l=0,\dots,m-1$, if $S=\{l\}$, we note
simply $|l\rangle$ for $|S\rangle$.

We have:
\begin{equation*}
  B\vert l\rangle = \sum_{l'=0}^{l} c^{l-l'} b^{m-1+l'-l} \vert l'\rangle +
  \sum_{l'=l+1}^{m-1} b^{l'-l-1} c^{m+l-l'}\vert l'\rangle.
\end{equation*}

For $z\in\mathbb{C}^*$, we define the vector
\begin{equation*}
  |z\rangle = \sum_l z^l |l\rangle.
\end{equation*}

Let us compute the action of $B$ on $|z\rangle$, by exchanging the sums over $l$
and $l'$:
\begin{align*}
  B\vert z \rangle &=
  \sum_{l'=0}^{m-1} \left(
  \sum_{l=l'}^{m-1}
  b^{m+l'-1-l} c^{l-l'} z^l
  +
  \sum_{l=0}^{l'-1}
  b^{l'-1-l} c^{m-l'+l} z^l
  \right)\vert l' \rangle\\
  &=
  \sum_{l'=0}^{m-1}
  \left(
  b^{m+l'} c^{-l'} \frac{(cz/b)^m - (cz/b)^{l'}}{cz-b}
  +
  b^{l'} c^{m-l'} \frac{(cz/b)^{l'} -1}{cz-b}
  \right)\vert l'\rangle \\
  &=
  \sum_{l'=0}^{m-1}
  \frac{c^m -b^m}{cz-b} z^{l'} \vert l'\rangle
  +
  \frac{z^m-1}{cz-b}
  \sum_{l'=0}^{m-1}
  b^{l'} c^{m-l'}\vert l'\rangle
\end{align*}

The first term is exactly $\frac{c^m-b^m}{cz-b}|z\rangle$. If we choose $z$ to
be a $m$th root of unity, then the factor in front of the second sum vanishes.
Then $|z\rangle$ (which is also an eigenvector of the translation operator) is
an eigenvector of $B$, with eigenvalue $\frac{c^m-b^m}{cz-b}$.
The eigenvalues are all distinct for $b\neq c$, and
$z=z_r=\exp\bigl(\frac{2ir\pi}{m}\bigr)$, with
$r=0,\dots,m-1$

These give $m=\binom{m}{1}$ orthogonal eigenvectors, and the one associated to
the largest eigenvalue $\frac{c^m-b^m}{c-b}$ corresponds to $z=1$ and has all
its entries equal.

\subsection{The $p=2$ sector}

If $S=\{l_1,l_2\}$, with $0\leq l_1 < l_2 \leq m-1$, we write $|l_1,l_2\rangle$
instead of $|S\rangle$. Let us write explicitly the action of $B$ on such a
vector.

As for $p=1$, there are two cases to consider:
\begin{itemize}
  \item either $0\leq l'_1 \leq l_1 < l'_2 \leq l_2 \leq m-1$,
  \item or $0\leq l_1 < l'_1 \leq l_2 < l'_2 \leq m-1$
\end{itemize}

Therefore the action of $B$ on the vector $\vert l_1,l_2\rangle$ can be splitted
into two sums:
\begin{multline*}
  B\vert l_1,l_2\rangle =
  \sum_{l'_1=0}^{l_1}\sum_{l'_2=l_1+1}^{l_2}
  c^{l_1-l'_1} b^{l'_2-l_1-1} c^{l_2-l'_2} b^{m+l'_1-l_2-1}
  \vert l'_1,l'_2 \rangle
  +\\
  \sum_{l'_1=l_1+1}^{l_2}\sum_{l'_2 = l_2+1}^{m-1}
  b^{l'_1-l_1-1} c^{l_2-l'_1} b^{l'_2-l_2-1} c^{m+l_1-l'_2}
  \vert l'_1,l'_2\rangle.
\end{multline*}

Define for $z_1, z_2\in\mathbb{C}^*$, the vector
\begin{equation*}
  |z_1,z_2\rangle = \sum_{l_1, l_2} z_1^{l_1} z_2^{l_2} |l_1, l_2\rangle,
\end{equation*}
where the summation is over all the allowed positions $0\leq l_1 < l_2 \leq
m-1$.
We compute $B\vert z_1,z_2\rangle$ and exchange sums over $l_i$'s and $l'_i$'s.
\begin{multline*}
  B\vert z_1, z_2\rangle =
  \sum_{l'_1<l'_2}
  \left[
    \left(
    \sum_{l_1=l'_1}^{l'_2-1}
    \sum_{l_2=l'_2}^{m-1}
    b^{m+l'_1+l'_2-l_1-l_2-2}
    c^{l_1-l'_1+l_2-l'_2}
    z_1^{l_1}z_2^{l_2}
    \right)
    \right.
    \\
    +
    \left.
    \left(
    \sum_{l_1=0}^{l'_1-1}
    \sum_{l_2=l'_1}^{l'_2-1}
    b^{l'_1+l'_2-l_1-l_2-2} c^{L+l_1+l_2-l'_1-l'_2}
    z_1^{l_1}z_2^{l_2}
    \right)
  \right]
  \vert l'_1, l'_2\rangle.
\end{multline*}
We compute explicitly the geometric series and obtain:
  \begin{multline*}
    B\vert z_1, z_2\rangle =
    \sum_{l'_1<l'_2}
    \left[
      b^{m+l'_1+l'_2}c^{-l'_1-l'_2}
      \frac{(cz_1/b)^{l'_2} - (cz_1/b)^{l'_1}}{c z_1 -b}
      \frac{(cz_2/b)^{m} - (cz_2/b)^{l'_2}}{c z_2 -b}
      +\right. \\
      \left.
      b^{l'_1+l'_2} c^{m-l'_1-l'_2}
      \frac{(cz_1/b)^{l'_1} - 1}{c z_1 -b}
      \frac{(cz_2/b)^{l'_2} - (cz_2/b)^{l'_1}}{c z_2 -b}
    \right]
  \vert l'_1, l'_2\rangle
  \end{multline*}

  We can factor the denominator $(c z_1 -b)(c z_2 -b)$. When we expand the
  numerators, we get eight terms which can be split into several categories:
  \begin{description}
  \item[wanted terms] those with $z_1^{l'_1} z_2^{l'_2}$:
    \begin{equation*}
      (b^m+c^{m})z_1^{l'_1} z_2^{l'_2},
    \end{equation*}
  \item[boundary terms] those appearing with a $z_i^{m}$ or $z_i^{0}=1$:
    \begin{equation*}
      -z_1^{l'_1} z_2^{m} b^{l'_2} c^{m-l'_2} - z_1^{0} z_2^{l'_2} b^{l'_1}
      c^{m-l'_1}
     + z_1^{l'_2} z_2^m b^{l'_1} c^{m-l'_1} + z_2^{l'_1} b^{l'_2} c^{m-l'_2},
    \end{equation*}
  \item[crossed terms] those appearing with $z_1$ and $z_2$ with the same power
    $l'_i$ ($i=1,2$):
    \begin{equation*}
      -z_1^{l'_2} z_2^{l'_2} b^{m+l'_1-l'_2} c^{l'_2-l'_1} - z_1^{l'_1} z_2^{l'_1}
      b^{l'_2-l'_1} c^{m+l'_1-l'_2}.
    \end{equation*}
  \end{description}

  We would like to get rid of terms of the last two categories, but keep those
  of the first one, which once resummed over $l_1'$ and $l_2'$ would look like
  the action of $B$ on an eigenvector.

  If instead of $\vert z_1,z_2\rangle$, we look at the antisymmetric version of it
  \begin{equation*}
    \vert \widetilde{z_1,z_2}\rangle = \vert z_1,z_2\rangle - \vert
    z_2,z_1\rangle,
  \end{equation*}
  the contribution of crossed terms will cancel, by symmetry considerations.

  When subtracting the antisymmetric counterpart (where $z_1$ and $z_2$ are
  exchanged) to boundary terms, the total
  contribution of boundary terms would cancel exactly,
  at the (sufficient) condition that $z_1$ and $z_2$ are $m$th root of $-1$ (and
  not 1 this time). There are $m$ such roots, but taking $z_1=z_2$ gives
  identically $0$ and
  exchanging $z_1$ and $z_2$ gives the same vector, up to a global sign. So
  there is $\binom{m}{2}$ choices, yielding distinct eigenvalues. The vectors
  $\vert \widetilde{z_1,z_2}\rangle$ form a basis of eigenvectors for the $p=2$
  sector.

  In this $p=2$ sector, the largest (positive) eigenvalue is
  $\frac{b^m+c^m}{|c\exp(i\pi/m)-b|^2}$,
  obtained by taking the two
  roots of $-1$ the closest to 1: $\{z_1, z_2\} = \{ \exp(i\pi/m),
  \exp(-i\pi/m)\}$.

\subsection{Sectors for general $p$}

For general $p$, we look for eigenvectors\footnote{%
  In more general form of the Bethe Ansatz, the coefficients $(-1)^\sigma$ are
  replaced by amplitude $C_\sigma$ depending in a more involved way on the
  permutation $\sigma$.
}
of the form:
\begin{equation*}
  \vert\widetilde{ z_1,\ldots z_{p}}\rangle =
  \sum_{l_1< \cdots < l_{p}}\sum_{\sigma\in\mathfrak{S}_p} (-1)^\sigma
    \prod_{j=1}^{p} z_{\sigma(j)} ^{l_j} |l_1,\ldots, l_{p}\rangle
  \end{equation*}
  where $\mathfrak{S}_p$ is the symmetric group over $p$ elements, and
  $z_1,\ldots,z_{p}\in\mathbb{C}^*$.

  When looking at the action of $B$ on this vector, then there are still terms
  of different types. The numeric factor in front of the wanted terms will have
  the form
  \begin{equation}
    \frac{c^m+(-1)^p b^m}{\prod_{j=0}^{p-1} (cz_j-b)}.
    \label{eq:eigenval_sect_p}
  \end{equation}
  The crossed terms (where two $z_j$ appear with the same exponent) cancel
  by antisymmetry, and the boundary terms will cancel if all the $z_j$ are $m$th
  roots of $(-1)^{p+1}$.

  Therefore the $\binom{m}{p}$ eigenvectors are obtained by choosing $p$
  distinct roots among the $m$th roots of $(-1)^{p+1}$.

  \begin{rem}
    It is important to notice that the
    eigenvectors do not depend on $b$ and $c$. In particular, they
    are also the eigenvectors for the transfer matrix $A$.
    Moreover, for generic values of $b$ and
    $c$, they are associated to disctinct eigenvalues: they are linearly
    independent and thus form a basis of the corresponding sector.
  \end{rem}

  Note that the largest eigenvalue of this sector, as in the case for
  $p=1$ and $p=2$, is obtained by
  choosing the $z_j$ to be the closest as possible to 1.

  When $p=2q$ is even, we chose the $z_j$ to be
  $\exp(\pm i\pi\frac{2r+1}{m})$, $r=0,\dots q-1$, and the associate eigenvalue is
  \begin{equation*}
    (c^m+b^m) \prod_{r=0}^{q-1} |c e^{i\pi\frac{2r+1}{m}} -b|^{-2}.
  \end{equation*}

  When $p=2q+1$ is odd, we chose the $z_j$ to be
  1, and $\exp(\pm i\pi\frac{2r}{m})$, $r=1,\dots q$, and the associate eigenvalue is
  \begin{equation*}
    \frac{c^m-b^m}{(c-b)}
    \prod_{r=1}^{q} |c e^{i\pi\frac{2r}{m}} -b|^{-2}.
  \end{equation*}

  \begin{rem}
    \label{rem:roots_unity}
    Denote by $U_{p,m}$ the set of $m$th roots of $(-1)^{p+1}$. Because of the
    fundamental identity satisfied by roots of unity:
    \begin{equation*}
      \prod_{z\in U_{p,m}} (b - cz) = b^m +(-1)^p c^m,
    \end{equation*}
    Equation~\eqref{eq:eigenval_sect_p} can be rewritten as
    \begin{equation*}
      \frac{c^m+(-1)^p}{\prod_{j=0}^{p-1}(c z_j -b)}
      =
      \prod_{\substack{z\in U_{p,m} \\ z\neq z_0,\dots,z_{p-1}} } (b-c z_j).
    \end{equation*}
  \end{rem}
\subsection{Taking the limit $b,c\to 1$}

By a continuity argument, the eigenvectors of the transfer matrix $A$ are those
computed above, and the associated eigenvalues are obtained by taking the limit
as $b$ and $c$ go to $1$ in \eqref{eq:eigenval_sect_p}.
In particular, the highest eigenvalue of the matrix $A$ in sector $p$ is given
by:
\begin{equation}
  \lambda_\text{max}^{(p)}=
  \begin{cases}
    \displaystyle
    2\prod_{r=0}^{q-1}
     |e^{i\pi\frac{2r+1}{2m}}
     2i\sin(\pi\frac{2r+1}{2m})|^{-2}
     = \frac{2}{\prod_{r=0}^{q-1} 4\sin^2(\pi\frac{2r+1}{2m})}
    & \text{if $p=2q$ is even,}\\
    \displaystyle
    m \prod_{r=1}^q
    |e^{i\pi\frac{r}{m}} 2i\sin(\pi\frac{r}{m})|^{-2}
    = \frac{m}{\prod_{r=1}^q 4 \sin^2(\frac{\pi r}{m})}.
 & \text{if $p=2q+1$ is odd.}
  \end{cases}
  \label{eq:lambda_max}
\end{equation}

Now notice that because of parity constraints, the number of vertices matched
with edges of one of the $m$-gons at the boundary is even. Thus the only values
of $p$ that one can see for $m$-barrel fullerenes are exactly those with the
same parity as $m$.

The leading contribution to the number of perfect matchings of $\Phi(F(m,k))$ as
$k$ goes to infinity is given,
up to corrections coming from the scalar product between the corresponding
normalized eigenvector and $|\Omega\rangle$,
by the $(k+1)$th power of the largest
of the largest eigenvalues of sectors of $A$ with $p$ and $m$ of same parity.

Since the function $\sin$ is increasing from $0$ to $\pi/2$, the value of
$\lambda_{\text{max}}^{(p)}$ is (for a given parity of $p$) is nondecreasing as a
function of $p$ as long as $\sin\bigl(\frac{\pi (p-1)}{2m}\bigr)$ is less or equal to
$\frac{1}{2}$, i.e. if $\frac{p-1}{2m}$ is less or equal to
$\frac{1}{6}$, As a consequence, the largest value $\lambda_{\text{max}}$ of the transfer
matrix $A$ is $\lambda_{\text{max}}^{(p_0)}$, for $p_0$ given by:
\begin{equation*}
  p_0=
  \left.
  \begin{cases}
    2\lfloor\frac{m+3}{6}\rfloor & \text{if $m$ is even} \\
    2\lfloor\frac{m}{6}\rfloor+1 & \text{if $m$ is odd}
  \end{cases}
\right\}
= m-2\left\lfloor \frac{m+1}{3}\right\rfloor.
\end{equation*}

One can indeed check that the leading terms of $\Phi(F(3,k))$, $\Phi(F(5,k))$
correspond to $p_0=1$ and the one for $\Phi(F(4,k))$ corresponds to $p_0=2$.
We are garanteed that the prefactor term coming from the scalar product is
nonzero since $|\Omega\rangle$ has nonnegative coefficients in the sector $p_0$
and the coefficients of the  eigenvector associated to $\lambda_{\text{max}}^{(p_0)}$
are all strictly positive.
So we can conclude that:

\begin{theo}
  \label{thm:result}
  For any fixed $m\geq 3$, the growth constant Equation~\eqref{eq:growth}
  for the family of $m$-barrel fullerenes $(F(m,k))_k$ is given by
  \begin{equation*}
    \rho(m)=
    \lambda_{\text{max}}^{(p_0)}.
  \end{equation*}
  with $p_0=m-2\left\lfloor
  \frac{m+1}{3}\right\rfloor$.
  Moreover, the number of horizontal rhombi on every slice converges in
  probability as $k$ goes to infinity to $p_0$.
\end{theo}

The proof of the last statement comes from the fact that for a fixed $m$,
the contribution of the sectors for $p\neq p_0$ to $\langle
\Omega|A^{k+1}|\Omega\rangle$ are exponentially small in $k$ when compared to
$\lambda_{\text{max}}^{(p_0)}$.

\section{Enumeration of perfect matchings of $F(m,k)$ through non-intersecting
paths}

Another alternative, now classical, technique to enumerate perfect matchings or
tilings, is
to encode the tiling with a family of non-intersecting paths, and use the
Lindstr\"om-Karlin-McGregor-Gessel-Viennot lemma~\cite{GesselViennot,KarlinMcGregor} to write the number of such paths as a
determinant which can then be evaluated. It turns out that the
family of paths corresponding to $m$-barrel fullerenes has been studied already
in great detail~\cite{Grabiner,KratBT}.

\begin{figure}[ht]
  \centering
  \begin{subfigure}[c]{4cm}
    \includegraphics[height=6cm]{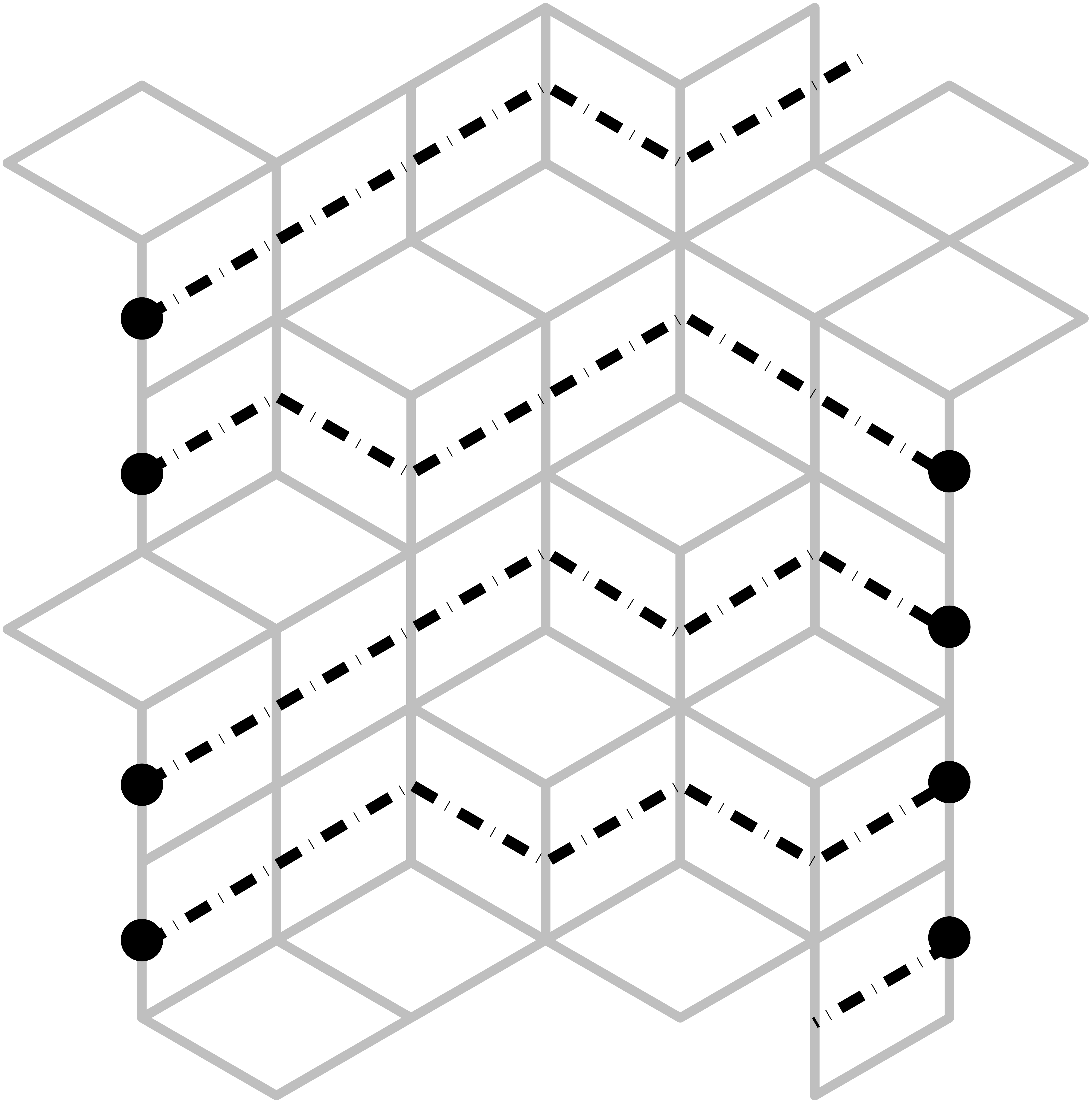}
  \end{subfigure}
  \hspace{0.2\textwidth}
  \begin{subfigure}[c]{4cm}
    \includegraphics[height=6cm]{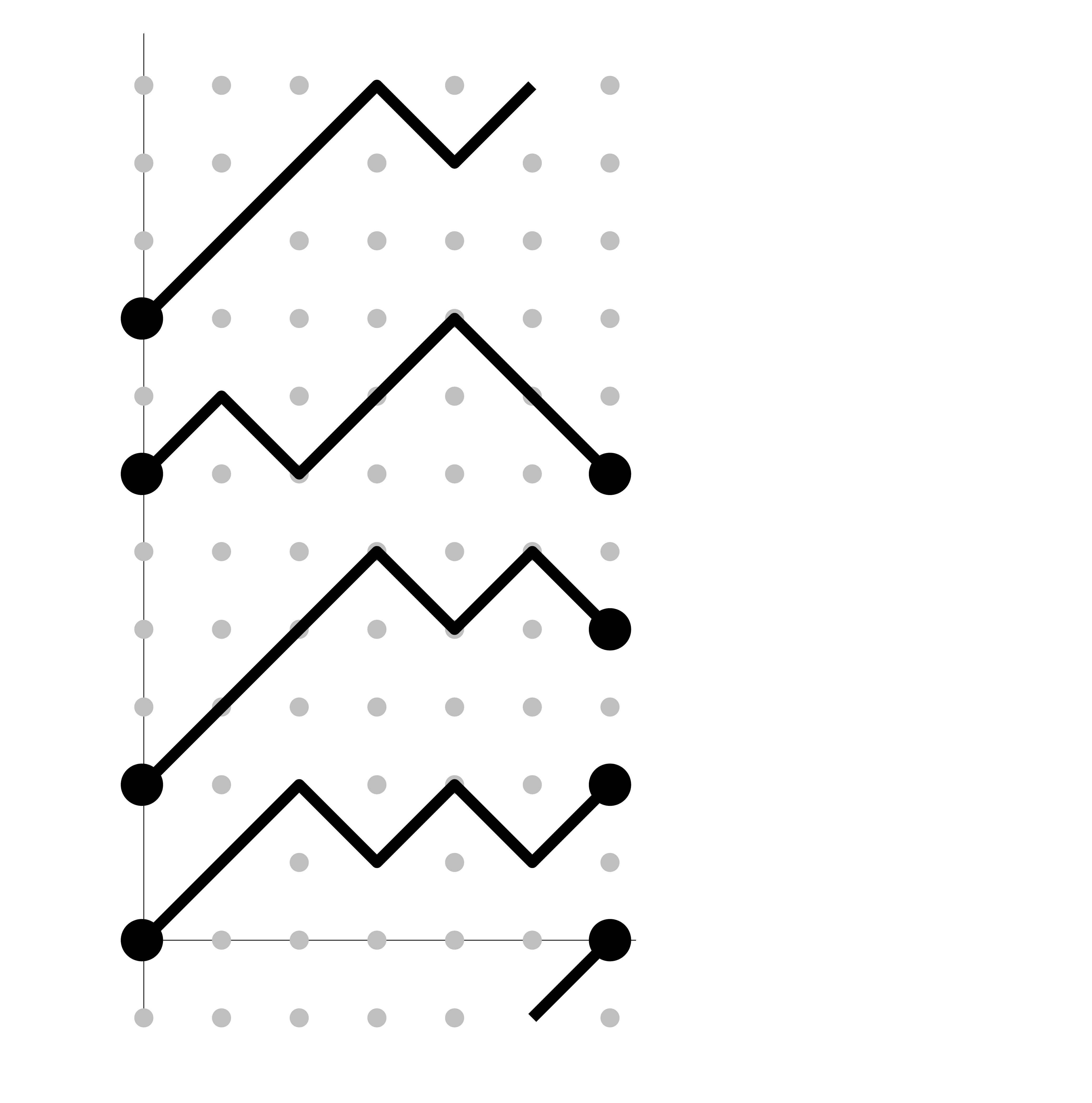}
  \end{subfigure}
\caption{Conversion of the rhombic tiling (left) to a family of non-intersecting
paths (right). Notice how the path starting from the left topmost starting point
winds around the cylinder to connect to the bottom most ending point on the
right.}
\label{barrel5}
\end{figure}

In order to relate our problem to existing results in the litterature, we need
to apply a last transformation to our rhombic tiling.
This is done in the following way: in the centers
of the vertical edges of the leftmost and rightmost non horizontal
rhombi, one places vertices, which will be the starting and ending points of the
paths. See Figure~\ref{barrel5}, where there are four vertices on each side.
Then one starts a path in each of the vertices on the left by connecting
midpoints of opposite vertical edges of rhombi. These paths necessarily
terminate in the midpoints of the vertical edges on the right.
In the Figure~\ref{barrel5}  on the left
these paths are marked by dotted lines. On the right, we have
slightly stretched vertically these paths so that they consist of up-steps
$(1,1)$ and down-steps $(1,-1)$. They have the property
that two paths do not have any point in common, and thus
form a family of non-intersecting paths viewed on the cylinder.

These paths on the right of Figure~\ref{barrel5}, can be understood as the
trajectories of a particle system evolving with time, flowing along the
horizontal axis. Each vertical line $x=t$ represents the particle configuration.
From time $t$ to time $t+1$, each of these particles
jumps one unit down or up. Initially, the mutual distances between particles are
even, and they will stay even at all times, but two particles will never be at
the same site at the same time.
Michael Fisher~\cite{Fisher} coined the term ``vicious walkers"
for this model of particles. More precisely, in his
model, we start with a fixed number of particles,
and at each time step, every particle moves one
unit in the positive or the negative direction so that at no time
two particles sit in the same place.
Since we started with a cylindrical graph, these particles do not actually move along
a line but actually along a circle with $2m$ sites.
Due to the particular matching problem that we started with,
namely, due to the fact that to the left and to the right of the
cylinder there is a ring of pentagons, the starting points of
the particles must come in pairs, with the particles in each
pair at distance~$2$, and the same applies to the end points.

Exact determinantal formulas have been given by Grabiner~\cite{Grabiner}
to count various families of such paths and Krattenthaler~\cite{KratBT} answered the
question of computing the asymptotics for this number of paths, as the number of
steps goes to infinity.
To stick with the notation of~\cite{KratBT}, positions of particles on the
circles will be labeled by half-integers. At $t=0$, and all the even times, the
positions of all particles will be (distinct) integers between 0 and $m-1$, and
at odd times, they will be odd multiples of $\frac{1}{2}$. If there are $k$
layers of hexagons, then the particle systems will evolve until time $k+1$.

Theorem~19 from~\cite{KratBT} give the number of families of $n$ paths with a
given starting and ending positions (or equivalently, the number of perfect
matchings of $F(m,k)$, with a fixed configuration on the two $m$-gons at
extremities). The starting and
ending positions are recorded in two sequences of numbers
$\eta=(\eta_1,\dots,\eta_n)$ and $\lambda=(\lambda_1,\dots,\lambda_n)$
respectively. We adapt it slightly to the context of $k+1$ time steps:

\begin{theo}[\cite{KratBT}, Theorem 19]
  Let $\eta=(\eta_1,\dots,\eta_n)$ be a vector of integers of half-integers with
  $m>\eta_1>\cdots > \eta_n\geq 0$ and $\lambda=(\lambda_1,\dots,\lambda_n)$ be
  a vector of integers of half-integers with $m>\lambda_{s+1} > \cdots >
  \lambda_n > \lambda_1 > \cdots >\lambda_s \geq 0$ for some $s$. Then, as $k$
  tends to $\infty$, in such a way that for all $j$, $k+1\equiv 2\eta_j +2\lambda_j\
  [2]$, the number of perfect matchings of $F(m,k)$ with boundary configurations
  encoded by $\eta$ and $\lambda$ is asymptotically equal to
  \begin{equation}
    \frac{2^{n^2-n}}{n m^n}
    \left(
    2^n \prod_{j=1}^n \cos \frac{\pi(j-\frac{n+1}{2})}{m}
    \right)^{k+1}
    \prod_{1\leq h<t\leq n}
    \left(%
    \sin\frac{\pi(\eta_h - \eta_t)}{m} \cdot
    \left|
    \sin\frac{\pi(\lambda_h - \lambda_t)}{m}
    \right|
    \right).
    \label{eq:asympt_krat}
  \end{equation}
\end{theo}



For getting the exact asymptotics for our problem, one would
simply have to sum this formula over all possible starting
and ending positions, that is, over all $\lambda$'s and $\eta$'s
which have the property that their coordinates come in pairs,
the two coordinates in each pair differing by~$1$ (cyclically)
and sum over $n$, as the number of paths is not fixed, but can be any even
number up to $m$.
Nevertheless, for fixed~$m$, we are talking about
a finite sum.

Clearly, not all choices of starting and ending points will
contribute to the leading term of the asymptotics. Inspection of
\eqref{eq:asympt_krat} shows that the relevant term in the formula is

\begin{equation*}
  \bigg( 2^{n}\prod _{j=1} ^{n}\cos\frac {\pi (j-\frac{n+1} {2})} {m}\bigg)^{k+1},
\end{equation*}
as everything else doesn't depend on~$k$.
In order to find the leading order of the asymptotics,
The task is thus to
find the (even) $n$ which maximises
\begin{equation}
  2^{n}\prod _{j=1} ^{n}\cos\frac {\pi(j-\frac {n+1} {2})} {m}.
  \label{eq:main_contrib_grabiner}
\end{equation}
But after closer inspection, we notice (unsurprisingly!) that, by
Remark~\ref{rem:roots_unity}\footnote{%
  More precisely, when taking the limit when $b$ and $c$ go to 1.
}, the quantity in Equation~\eqref{eq:main_contrib_grabiner} is equal
to $\lambda_\text{max}^{(p)}$ from \eqref{eq:lambda_max} with $n+p=m$, and the
value of $n$ maximizing~\eqref{eq:main_contrib_grabiner} is
$n_0=m-p_0=2\lfloor \frac{m+1}{3}\rfloor$.

%
%


We thus recovered from another method the value of the growth constant
\begin{equation*}
  \rho(m)=2^{2\lfloor\frac{m+1}{3}\rfloor}
  \prod _{j=1} ^{2\lfloor\frac{m+1}{3}\rfloor}\cos\frac {\pi(j-\frac {n+1} {2})} {m}
  =
  \prod _{j=1} ^{\lfloor(m+1)/3\rfloor}
  \left(2\cos\frac {\pi(2j-1)} {2m}\right)^2,
\end{equation*}
as in Theorem~\ref{thm:result}.


 \section{The entropy of the family $(F(m,k))$}

Recall that the graph $F(m,k)$ has $2m(k+2)$ vertices. The \emph{dimer entropy}
~\cite{FKLM08, FP05}
$h(m)$ of the family $(F(m,k))_{k\geq 0}$ is defined as
\begin{equation*}
  h(m):=\limsup_{k\to\infty} \frac{\log\Phi (F(m,k))}{2m(k+2)}.
\end{equation*}
 From the previous sections we deduce
 \begin{equation*}
   h(m)=\frac{\log \rho(m)}{2m}.
 \end{equation*}
 Equivalently, it says that the number of matchings if $F(m,k)$ for $m$ fixed
and $k\gg 1$ is of order $e^{kh(m)}$.

The quantity $h(m)$ can be written as a Riemann sum, which increases and
converges to the corresponding integral as $m$ goes to infinity:
\begin{equation*}
  \mathbf{h}=\frac{-3}{2\pi} \int_0^{\pi/3} \log (2\sin t)\mathrm{d}t
  \simeq 0.1615329736\dots
\end{equation*}
which is the maximal entropy for ergodic Gibbs measure on uniform dimer
configurations of the infinite hexagonal lattice~\cite{Kenyon1,Kenyon2}.


As $m=5$ or $6$, we see that the entropy for the family $F(m,k)$ is strictly
smaller than $\mathbf{h}$. This is due to the geographic localisation of the
pentagons in our family of fullerenes. However,
if the pentagons are far enough from one another, their presence
pentagons seems very unlikely to have a great
nfluence on the number of perfect matchings of a large fullerene which would
behave in this respect just as a the hexagonal lattice.
That motivates us to advance the following conjecture:

Let $n$ be an integer number greater than $11$ so that there is a fullerene
with $2n$ vertices. Denote by $\mu_{2n}$ the maximal number of perfect
matchings in all fullerene graphs with $2n$ vertices.  Define
 \[h_F:=\limsup_{n\to\infty} \frac{\log\mu_{2n}}{2n}.\]
 We conjecture that $h_F=\mathbf{h}$.

A similar claim seems plausible also for $m$-generalized fullerenes. Fix an
integer $m\ge 3, m\ne 5,6$. Let $\mu_{2m(k+2),m}$ be the
maximal number of perfect matchings
 in all $m$-generalized fullerene graphs with $2m(k+2)$ vertices.  Define
 \[h_F(m):=\limsup_{k\to\infty} \frac{\log\mu_{2m(k+2),m}}{2m(k+2)}.\]
 We conjecture that $h_F(m)$ is also equal to $\mathbf{h}$ (and thus strictly
 greater than $h(m)$).

\section*{Acknowledgement}

This work has been supported by the Center for International Scientific Studies
and Collaboration (CISSC) and French Embassy in Tehran. We warmly thank
Christian Krattenthaler for very valuable discussions about this topic.

%

%
\end{document}

%% file: notation_edges.pdf_t
\begin{picture}(0,0)%
\includegraphics{notation_edges.pdf}%
\end{picture}%
\setlength{\unitlength}{4144sp}%
\begingroup\makeatletter\ifx\SetFigFont\undefined%
\gdef\SetFigFont#1#2#3#4#5{%
  \reset@font\fontsize{#1}{#2pt}%
  \fontfamily{#3}\fontseries{#4}\fontshape{#5}%
  \selectfont}%
\fi\endgroup%
\begin{picture}(5909,8188)(1336,-7755)
\put(1351,-511){\makebox(0,0)[lb]{\smash{{\SetFigFont{34}{40.8}{\familydefault}{\mddefault}{\updefault}{\color[rgb]{0,0,0}$e_{j-1,l+1}$}%
}}}}
\put(1351,-5911){\makebox(0,0)[lb]{\smash{{\SetFigFont{34}{40.8}{\familydefault}{\mddefault}{\updefault}{\color[rgb]{0,0,0}$e_{j-1,l}$}%
}}}}
\put(4951,-3211){\makebox(0,0)[lb]{\smash{{\SetFigFont{34}{40.8}{\familydefault}{\mddefault}{\updefault}{\color[rgb]{0,0,0}$e_{j,l}$}%
}}}}
\end{picture}%

%% file: asymptoticfullerene.bbl
\begin{thebibliography}{19}
\providecommand{\natexlab}[1]{#1}
\providecommand{\url}[1]{\texttt{#1}}
\expandafter\ifx\csname urlstyle\endcsname\relax
  \providecommand{\doi}[1]{doi: #1}\else
  \providecommand{\doi}{doi: \begingroup \urlstyle{rm}\Url}\fi

\bibitem[Andova et~al.(2012)Andova, Do{\v s}li\'c, Krnc, Lu{\v z}ar, and {\v
  S}krekovski]{andova1}
V.~Andova, T.~Do{\v s}li\'c, M.~Krnc, B.~Lu{\v z}ar, and R.~{\v S}krekovski.
\newblock On the diameter and some related invariants of fullerene graphs.
\newblock \emph{MATCH Commun. Math. Comput. Chem.}, 68\penalty0 (1):\penalty0
  109--130, 2012.
\newblock ISSN 0340-6253.

\bibitem[Behmaram and Friedland(2013)]{BF12}
A.~Behmaram and S.~Friedland.
\newblock Upper bounds for perfect matchings in {P}faffian and planar graphs.
\newblock \emph{Electron. J. Combin.}, 20\penalty0 (1):\penalty0 Paper 64, 16,
  2013.
\newblock ISSN 1077-8926.

\bibitem[Behmaram et~al.(2016)Behmaram, Do{\v s}li\'c, and Friedland]{BF1}
A.~Behmaram, T.~Do{\v s}li\'c, and S.~Friedland.
\newblock Matchings in {$m$}-generalized fullerene graphs.
\newblock \emph{Ars Math. Contemp.}, 11\penalty0 (2):\penalty0 301--313, 2016.
\newblock ISSN 1855-3966.

\bibitem[Bethe(1931)]{Bethe1931}
H.~Bethe.
\newblock Zur theorie der metalle.
\newblock \emph{Zeitschrift f{\"u}r Physik}, 71\penalty0 (3):\penalty0
  205--226, Mar 1931.
\newblock ISSN 0044-3328.
\newblock \doi{10.1007/BF01341708}.
\newblock URL \url{https://doi.org/10.1007/BF01341708}.

\bibitem[Do{\v{s}}li{\'{c}}(2008)]{dosliclff}
T.~Do{\v{s}}li{\'{c}}.
\newblock Leapfrog fullerenes have many perfect matchings.
\newblock \emph{Journal of Mathematical Chemistry}, 44\penalty0 (1):\penalty0
  1--4, Jul 2008.
\newblock ISSN 1572-8897.
\newblock \doi{10.1007/s10910-007-9287-x}.
\newblock URL \url{https://doi.org/10.1007/s10910-007-9287-x}.

\bibitem[Do{\v{s}}li{\'{c}}(2009)]{doslicfmm}
T.~Do{\v{s}}li{\'{c}}.
\newblock Finding more perfect matchings in leapfrog fullerenes.
\newblock \emph{Journal of Mathematical Chemistry}, 45\penalty0 (4):\penalty0
  1130--1136, Apr 2009.
\newblock ISSN 1572-8897.
\newblock \doi{10.1007/s10910-008-9435-y}.
\newblock URL \url{https://doi.org/10.1007/s10910-008-9435-y}.

\bibitem[Fisher(1984)]{Fisher}
M.~E. Fisher.
\newblock Walks, walls, wetting, and melting.
\newblock \emph{J. Statist. Phys.}, 34\penalty0 (5-6):\penalty0 667--729, 1984.
\newblock ISSN 0022-4715.
\newblock \doi{10.1007/BF01009436}.
\newblock URL \url{http://dx.doi.org/10.1007/BF01009436}.

\bibitem[Friedland and Peled(2005)]{FP05}
S.~Friedland and U.~N. Peled.
\newblock Theory of computation of multidimensional entropy with an application
  to the monomer{\textendash}dimer problem.
\newblock \emph{Advances in Applied Mathematics}, 34\penalty0 (3):\penalty0
  486--522, apr 2005.
\newblock \doi{10.1016/j.aam.2004.08.005}.
\newblock URL \url{https://doi.org/10.1016/j.aam.2004.08.005}.

\bibitem[Friedland et~al.(2008)Friedland, Krop, Lundow, and
  Markstr{\"o}m]{FKLM08}
S.~Friedland, E.~Krop, P.~H. Lundow, and K.~Markstr{\"o}m.
\newblock On the validations of the asymptotic matching conjectures.
\newblock \emph{Journal of Statistical Physics}, 133\penalty0 (3):\penalty0
  513--533, Nov 2008.
\newblock ISSN 1572-9613.
\newblock \doi{10.1007/s10955-008-9550-y}.
\newblock URL \url{https://doi.org/10.1007/s10955-008-9550-y}.

\bibitem[Gessel and Viennot(1985)]{GesselViennot}
I.~Gessel and G.~Viennot.
\newblock {Binomial determinants, paths, and hook length formulae}.
\newblock \emph{Adv. in Math.}, 58\penalty0 (3):\penalty0 300--321, 1985.
\newblock ISSN 0001-8708.

\bibitem[Grabiner(2002)]{Grabiner}
D.~J. Grabiner.
\newblock Random walk in an alcove of an affine {W}eyl group, and non-colliding
  random walks on an interval.
\newblock \emph{J. Combin. Theory Ser. A}, 97\penalty0 (2):\penalty0 285--306,
  2002.
\newblock ISSN 0097-3165.
\newblock \doi{10.1006/jcta.2001.3216}.
\newblock URL \url{http://dx.doi.org/10.1006/jcta.2001.3216}.

\bibitem[Gr\"unbaum and Motzkin(1963)]{grunbaum}
B.~Gr\"unbaum and T.~S. Motzkin.
\newblock The number of hexagons and the simplicity of geodesics on certain
  polyhedra.
\newblock \emph{Canad. J. Math.}, 15:\penalty0 744--751, 1963.
\newblock ISSN 0008-414X.
\newblock \doi{10.4153/CJM-1963-071-3}.
\newblock URL \url{http://dx.doi.org/10.4153/CJM-1963-071-3}.

\bibitem[Kardo{\v{s}} et~al.(2009)Kardo{\v{s}}, Kr{\'a}l', Mi{\v{s}}kuf, and
  Sereni]{KKMS09}
F.~Kardo{\v{s}}, D.~Kr{\'a}l', J.~Mi{\v{s}}kuf, and J.-S. Sereni.
\newblock Fullerene graphs have exponentially many perfect matchings.
\newblock \emph{Journal of Mathematical Chemistry}, 46\penalty0 (2):\penalty0
  443--447, Aug 2009.
\newblock ISSN 1572-8897.
\newblock \doi{10.1007/s10910-008-9471-7}.
\newblock URL \url{https://doi.org/10.1007/s10910-008-9471-7}.

\bibitem[Karlin and McGregor(1959)]{KarlinMcGregor}
S.~Karlin and J.~McGregor.
\newblock {Coincidence probabilities}.
\newblock \emph{Pacific J. Math.}, 9:\penalty0 1141--1164, 1959.
\newblock ISSN 0030-8730.

\bibitem[Kenyon(1997)]{Kenyon1}
R.~Kenyon.
\newblock {Local statistics of lattice dimers}.
\newblock \emph{Ann. Inst. H. Poincar{\'e} Probab. Statist.}, 33\penalty0
  (5):\penalty0 591--618, 1997.
\newblock ISSN 0246-0203.

\bibitem[Kenyon(2002)]{Kenyon2}
R.~Kenyon.
\newblock {The {L}aplacian and {D}irac operators on critical planar graphs}.
\newblock \emph{Invent. Math.}, 150\penalty0 (2):\penalty0 409--439, 2002.
\newblock ISSN 0020-9910.

\bibitem[Krattenthaler(2004/07)]{KratBT}
C.~Krattenthaler.
\newblock Asymptotics for random walks in alcoves of affine {W}eyl groups.
\newblock \emph{S\'em. Lothar. Combin.}, 52:\penalty0 Art. B52i, 72, 2004/07.
\newblock ISSN 1286-4889.

\bibitem[{Kroto} et~al.(1985){Kroto}, {Heath}, {Obrien}, {Curl}, and
  {Smalley}]{kroto}
H.~W. {Kroto}, J.~R. {Heath}, S.~C. {Obrien}, R.~F. {Curl}, and R.~E.
  {Smalley}.
\newblock {C(60): Buckminsterfullerene}.
\newblock \emph{Nature}, 318:\penalty0 162, Nov. 1985.
\newblock \doi{10.1038/318162a0}.

\bibitem[Kutnar and Maru{\v s}i{\v c}(2008)]{KM08}
K.~Kutnar and D.~Maru{\v s}i{\v c}.
\newblock On cyclic edge-connectivity of fullerenes.
\newblock \emph{Discrete Appl. Math.}, 156\penalty0 (10):\penalty0 1661--1669,
  2008.
\newblock ISSN 0166-218X.
\newblock \doi{10.1016/j.dam.2007.08.046}.
\newblock URL \url{http://dx.doi.org/10.1016/j.dam.2007.08.046}.

\end{thebibliography}
